\documentclass[12pt,a4paper]{article}
\usepackage{natbib}
\usepackage{amsmath,amssymb,graphicx,time}
\usepackage[colorlinks,linkcolor=blue,anchorcolor=blue,citecolor=blue,urlcolor=blue,pdfauthor={Keith Briggs},pdfmenubar=false,pdffitwindow=true,pdfwindowui=false,pdftitle={Poisson_maxima}]{hyperref}

\def\pr{{\text{Pr}}}
\newcommand*{\Prob}[1]{\ensuremath{\text{Pr}[#1]}}

\newcommand{\di}{\text{d}}

\newcommand{\R}{\mathbb}
\def\R2{{$\mathbb{R}^2$}}

\def\r2{{\mathbb{R}^2}}

\title{\vspace{-10mm}A note on the distribution of the maximum of a set of Poisson random variables}
\author{Keith Briggs\thanks{corresponding author; \texttt{keith.briggs@bt.com}},  Linlin Song (BT Research, Martlesham)\\ \& Thomas Prellberg (Mathematics, QMUL)}

\def\today{\number\year\space
 \ifcase\month\or January\or February\or March\or April\or May\or June\or
   July\or August\or September\or October\or November\or December\fi
     \space\number\day
}
\date{\today\ \now}
\date{2009-03-12} 

\def\SourceFileFooter{{\noindent\tiny Typeset with pdf\LaTeXe. Source \url{gold:~kbriggs/Poisson_maxima-0.1/Poisson_maxima_03.tex} last changed 2009-03-12 11:19}}

\normalsize

\begin{document}
\maketitle

\vspace{-5mm}\begin{center}
\begin{minipage}{10cm}
Given a set of independent Poisson random variables%
with common mean,  we study the distribution of their maximum 
and obtain an accurate asymptotic formula to locate the most probable value of  the maximum. 
We verify our analytic results with very precise numerical computations.
\end{minipage}
\end{center}

\noindent We deal with a set of independent Poisson random variables
$\{X_1,X_2,\dots, X_n\}$ with common mean $\lambda$, so that
$\pr[X_{i}=k]={e^{-\lambda}\lambda^{k}}/{k!}$.  We let $M_{n}=\max(X_i)$ and
wish to describe the distribution of $M_n$.  Our motivation is a problem in
random graph theory, where we were interested in the distribution of maximum
degree in graphs with Poisson degree distribution.

\begin{figure}%
\begin{center}
\includegraphics[width=0.243\hsize]{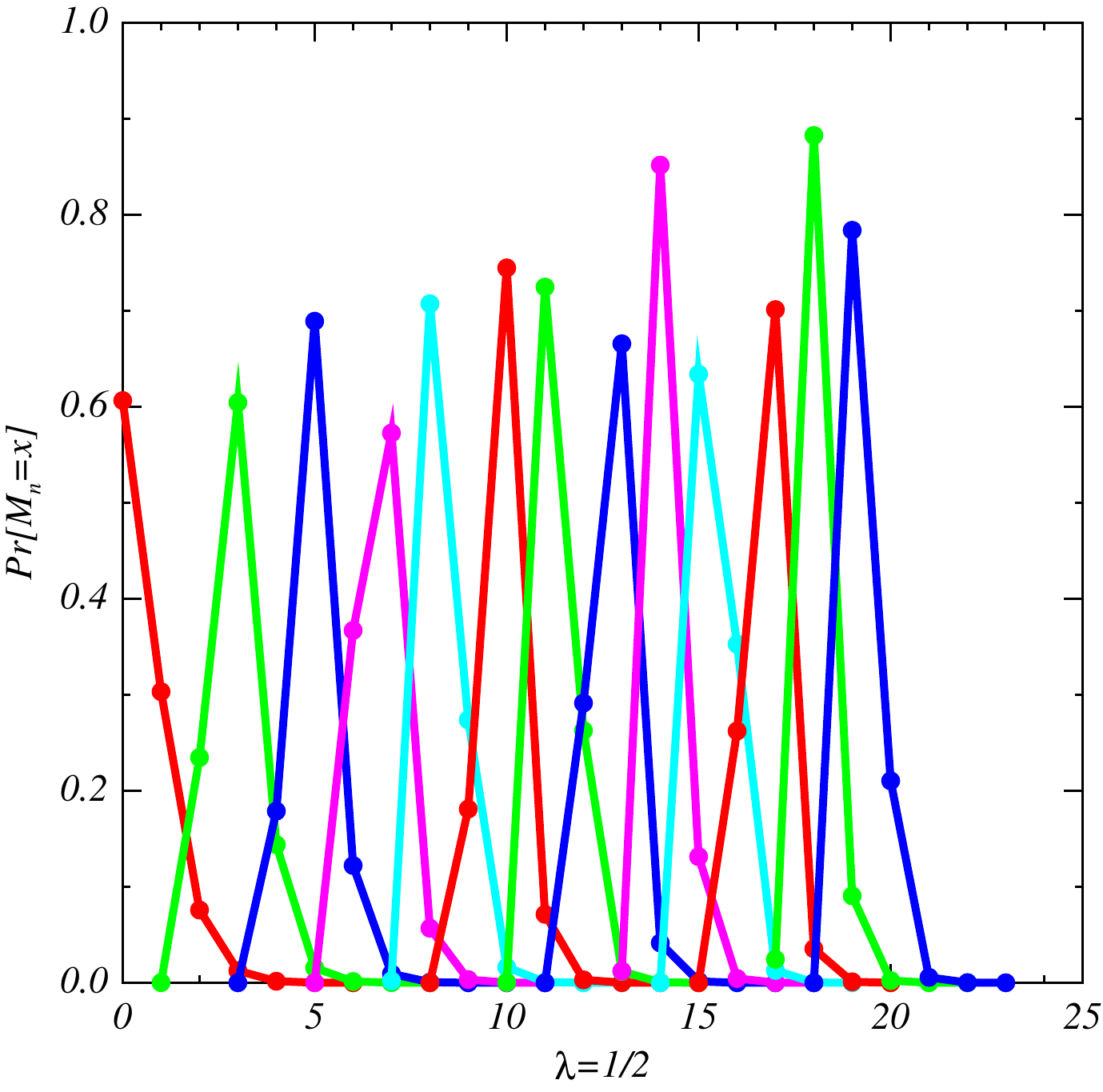}
\includegraphics[width=0.23\hsize]{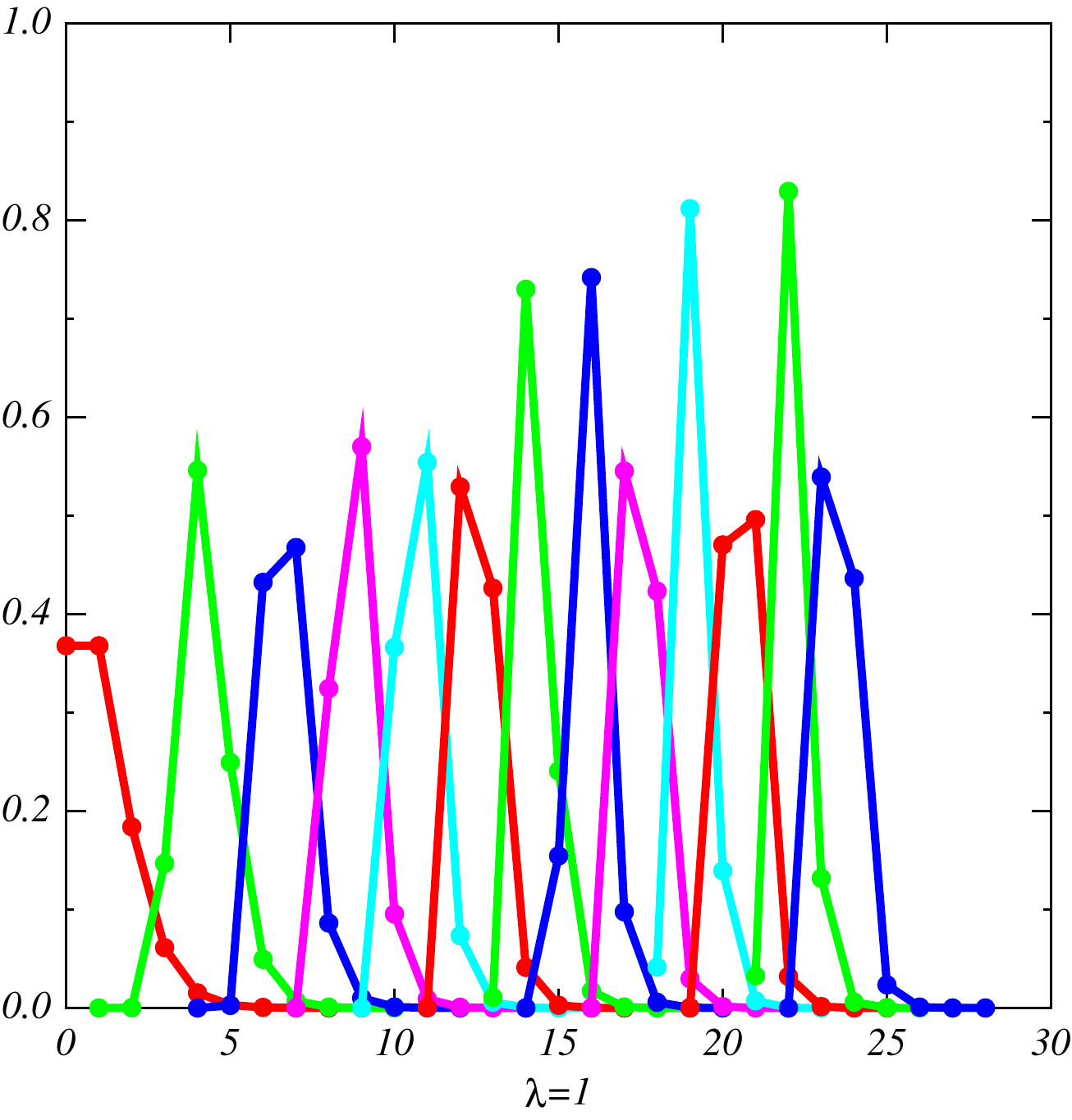}
\includegraphics[width=0.23\hsize]{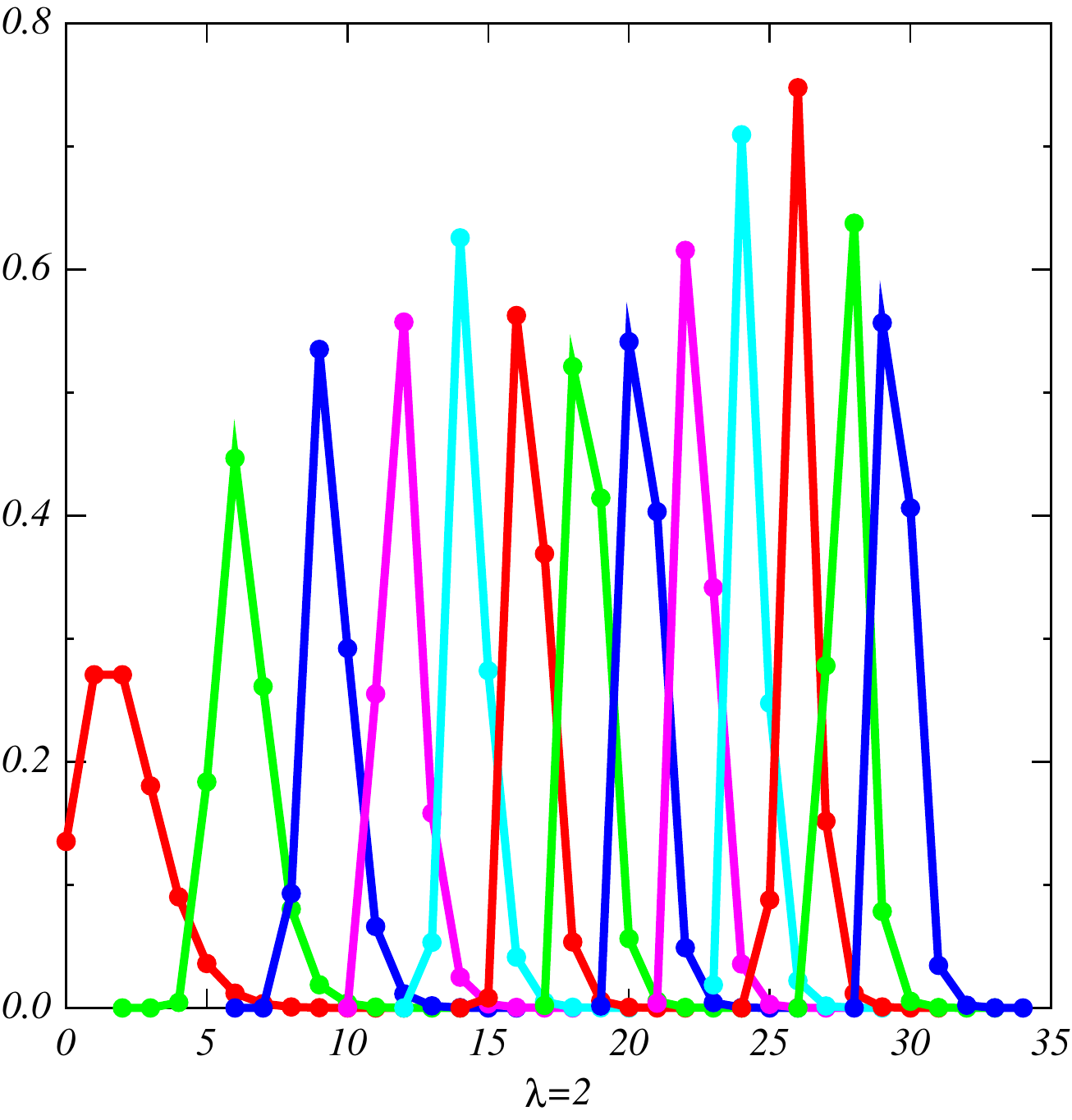}
\includegraphics[width=0.23\hsize]{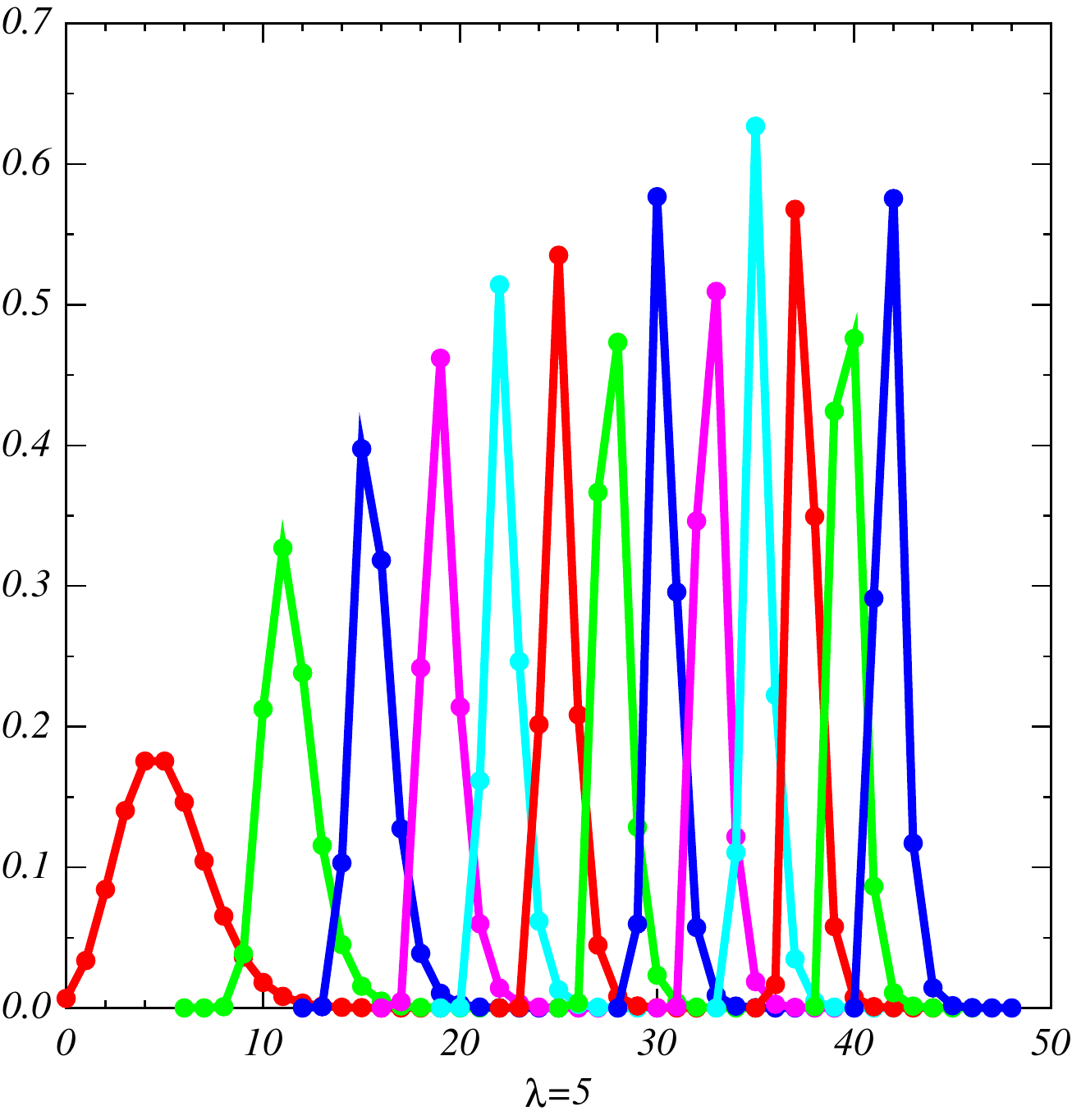}
\end{center}
\caption{From left to right: the distribution of the maximum of Poisson
variables for $\lambda=1/2, 1, 2, 5$ (left to right) and $n=10^0, 10^2, 10^4, \dots,
10^{24}.$ Note that there is an error in Fig~1 in \citet{Anderson1997}, where the curves labelled $k=6$ and $k=8$ are incorrect. 
}
\label{distribution}
\end{figure}

We have
\begin{eqnarray}
  \Prob{X_i<k}  = Q(k,\lambda)\equiv\Gamma(k,\lambda)/\Gamma(k)\nonumber
\end{eqnarray}
where $Q$ and $\Gamma(\cdot,\cdot)$ are incomplete Gamma functions; that is,
\begin{equation}
 \Gamma(a,x)=\int_{x}^{\infty}t^{a-1}e^{-t}\di t.\nonumber
\end{equation}
From the independence of the Poisson variables,
\begin{equation}\label{ccdf}
 \Prob{M_n\leqslant k}=\Prob{X_1 \leqslant k}^n=Q(k+1,\lambda)^n={\Gamma(k+1,\lambda)^n/\Gamma(k+1)}^n.\nonumber
\end{equation}
Our aim is to approximate the distribution of $M_n$.  We have
\begin{eqnarray}\label{PrMn=k}
  \Prob{M_n=k}&=&\Prob{M_n\leqslant k}-\Prob{M_n\leqslant
  k-1}\nonumber\\
  &=&Q(k+1,\lambda)^n-Q(k,\lambda)^n\nonumber
\end{eqnarray}
Examples of these distributions are shown in Figure~\ref{distribution}. These
numerical results demonstrate the so-called \textit{focussing} effect; the
maxima $M_n$ are concentrated on at most two adjacent integers for large $n$;
we call them \textit{modal values}.  It is this focussing that allows us to
characterize the distributions very precisely by a single asymptotic estimate.

\begin{figure}%
\begin{center}
\includegraphics[width=0.243\hsize]{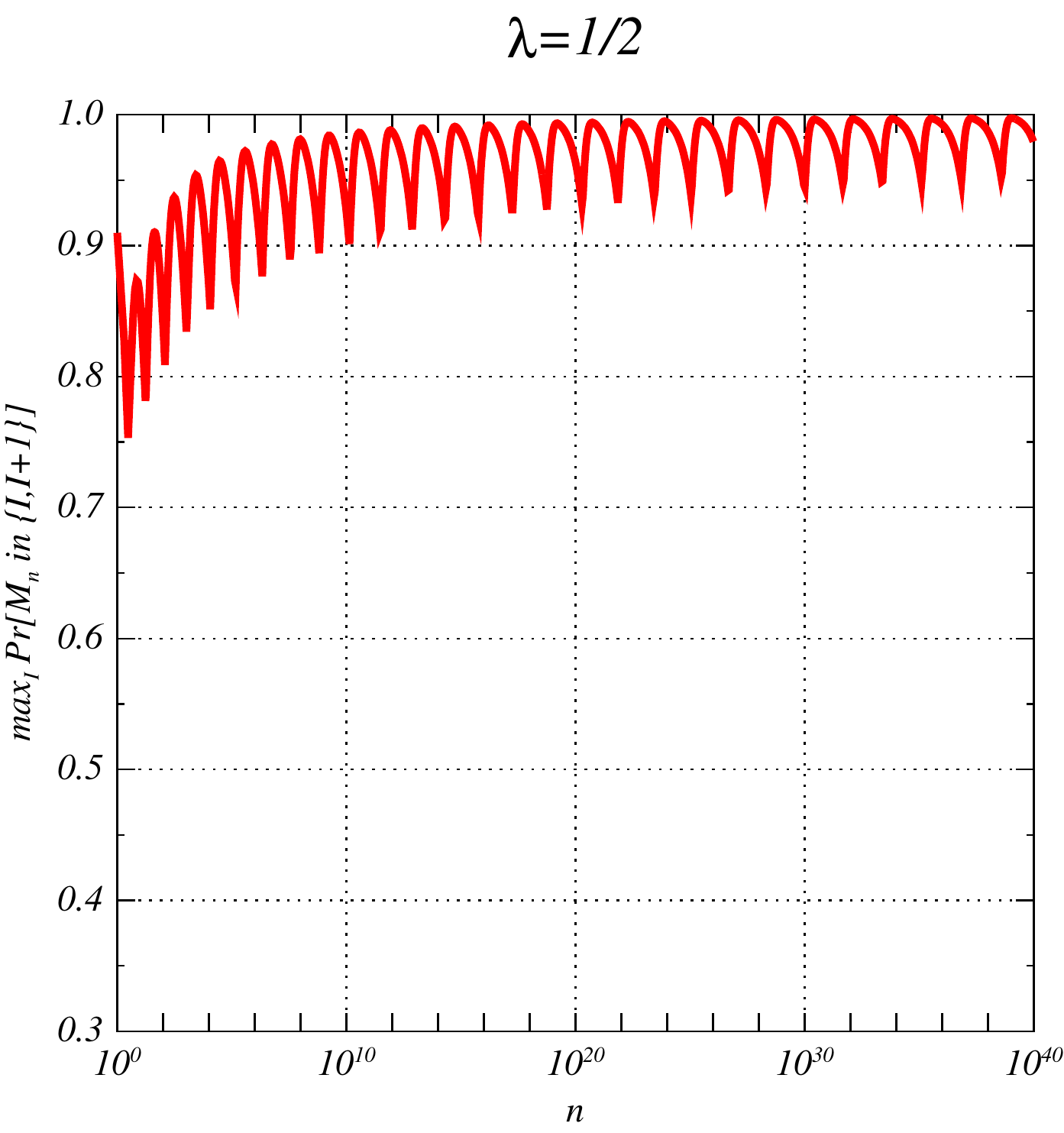}
\includegraphics[width=0.23\hsize]{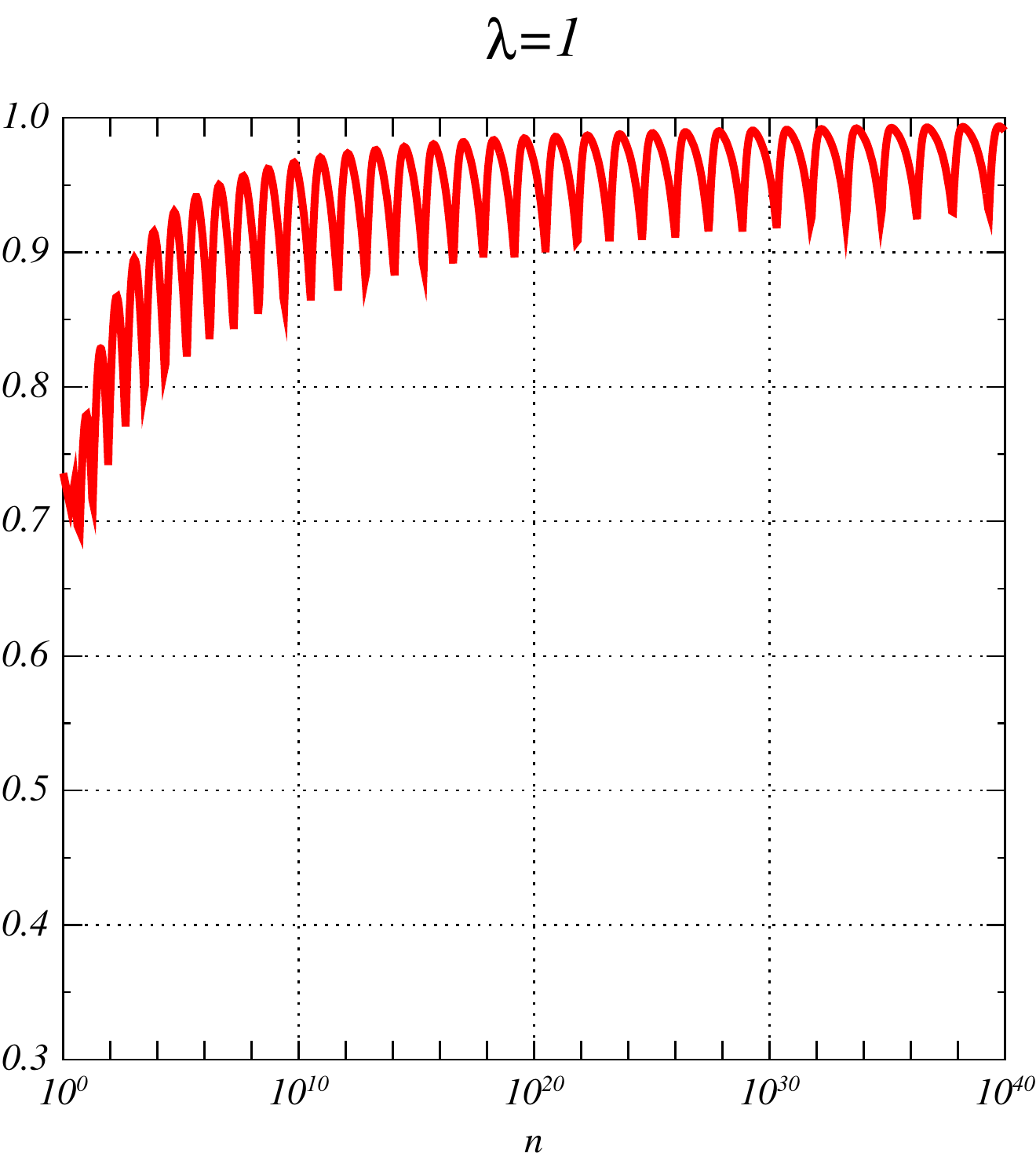}
\includegraphics[width=0.23\hsize]{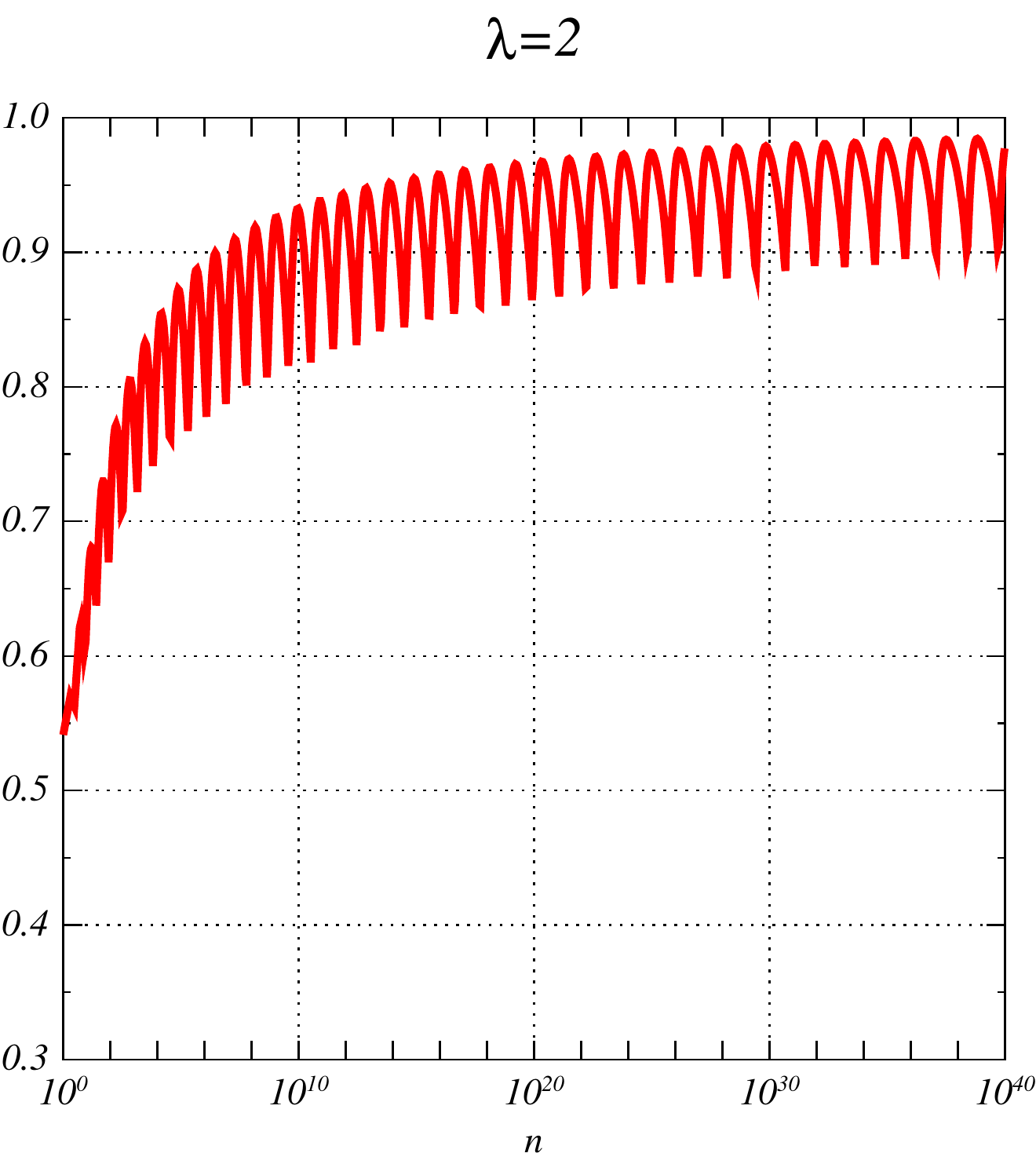}
\includegraphics[width=0.23\hsize]{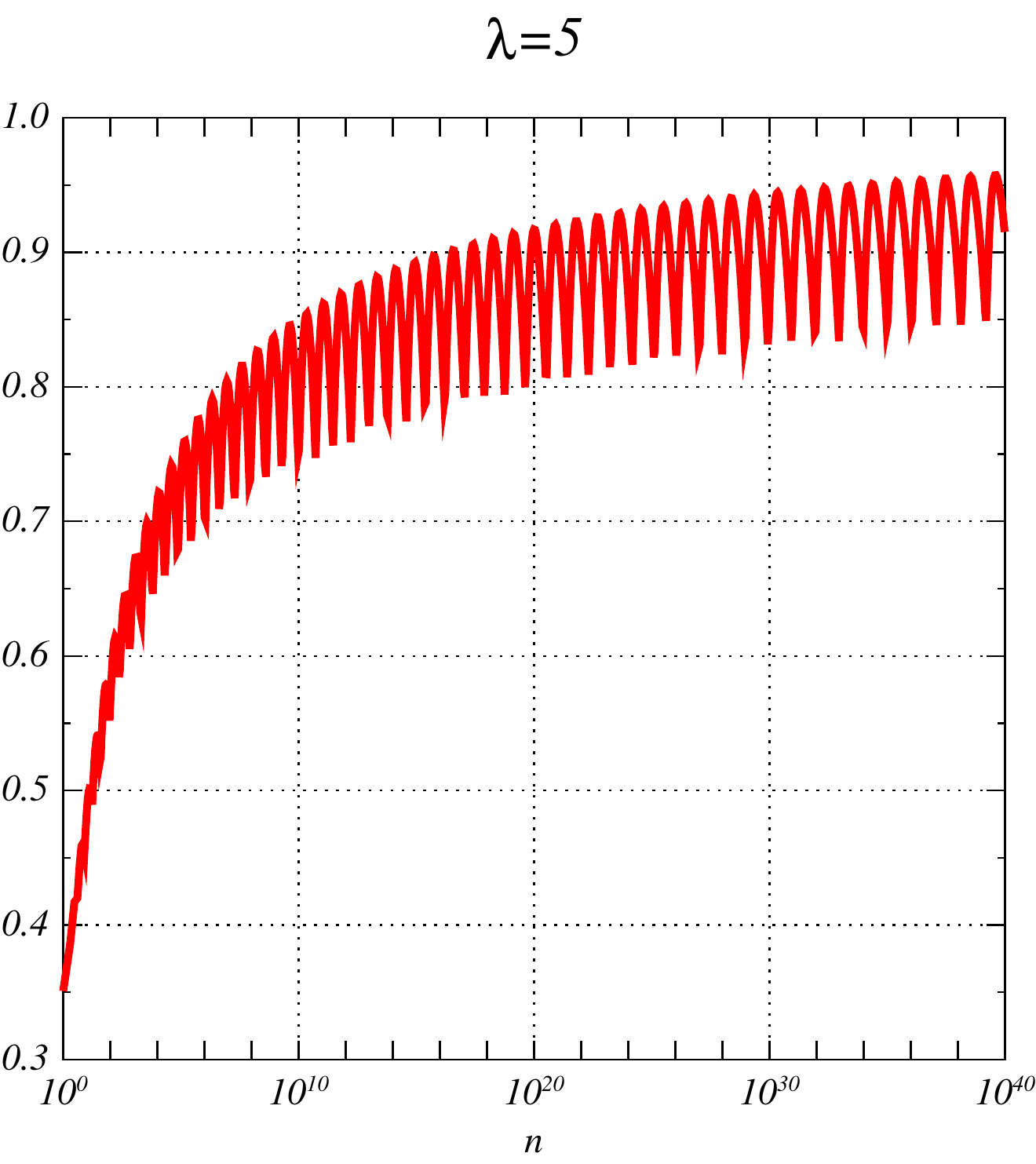}
\end{center}
\caption{The maximal probability (with respect to $I_n$) that $M_n \in \{I_n,I_n+1\}$ for $\lambda=1/2, 1, 2, 5$  (left to right)  and $10^0\leqslant n\leqslant 10^{40}$.  The curves show the probability that $M_n$ takes either of its two most frequently occurring values.}\label{prob}
\end{figure}

In previous work on this problem, \citet{Anderson1970} proved the existence of
integers $I_n$ such that $\Prob{M_n \in (I_n,I_n+1)} \to 1 $ as $n\rightarrow\infty$
for fixed $\lambda>0$; and that $I_n\sim\beta_n$, where $\beta_n$ is defined as
the unique solution of $Q(\beta_n,\lambda)=1/n$.

Following this work, \citet{Kimber1983} computed an asymptotic result;
he showed $I_n\sim\log n/{\log\log n}$ and $P_n \sim (k/{I_n})^{1+B_n}$ with
$B_n$ dense in $[-1/2, 1/2]$, and concluded that to the first order, the rate
of growth of $I_n$ is independent of the Poisson parameter $\lambda$. He
concluded that $P_n$, defined as $P_n=\Prob{M_n \in (I_n, I_n+1)}$,
oscillates and the oscillation persists for arbitrarily large $n$. We
illustrate in Figure~\ref{prob} exactly how this probability oscillates as $n
\to \infty$.  Our numerical experiments show that $\log n/\log\log n$
estimates $I_n$ very poorly. We aim to improve this asymptotic formula.

Our method is a refinement of that of Kimber; that is, we consider a continuous
distribution $g$ which interpolates the Poisson maximum distribution, and we
solve $g(x)=1/n$.   Consider $g_\lambda(x)\equiv
1-\Gamma(x+1,\lambda)/\Gamma(x+1)$ for fixed $\lambda \in \mathbb{R}^{+}$,
which is a strictly decreasing function on $(0,\infty)$.   If $\epsilon=1/n$ is
a small positive real, then $g_\lambda(x)$ has a unique root $x(\epsilon)$
which increases as $\epsilon \to 0^{+}$.  We will develop an asymptotic
expansion (as $\epsilon \to 0$) of this root $x(\epsilon)$.

\begin{figure}%
\begin{center}
\includegraphics[width=0.243\linewidth]{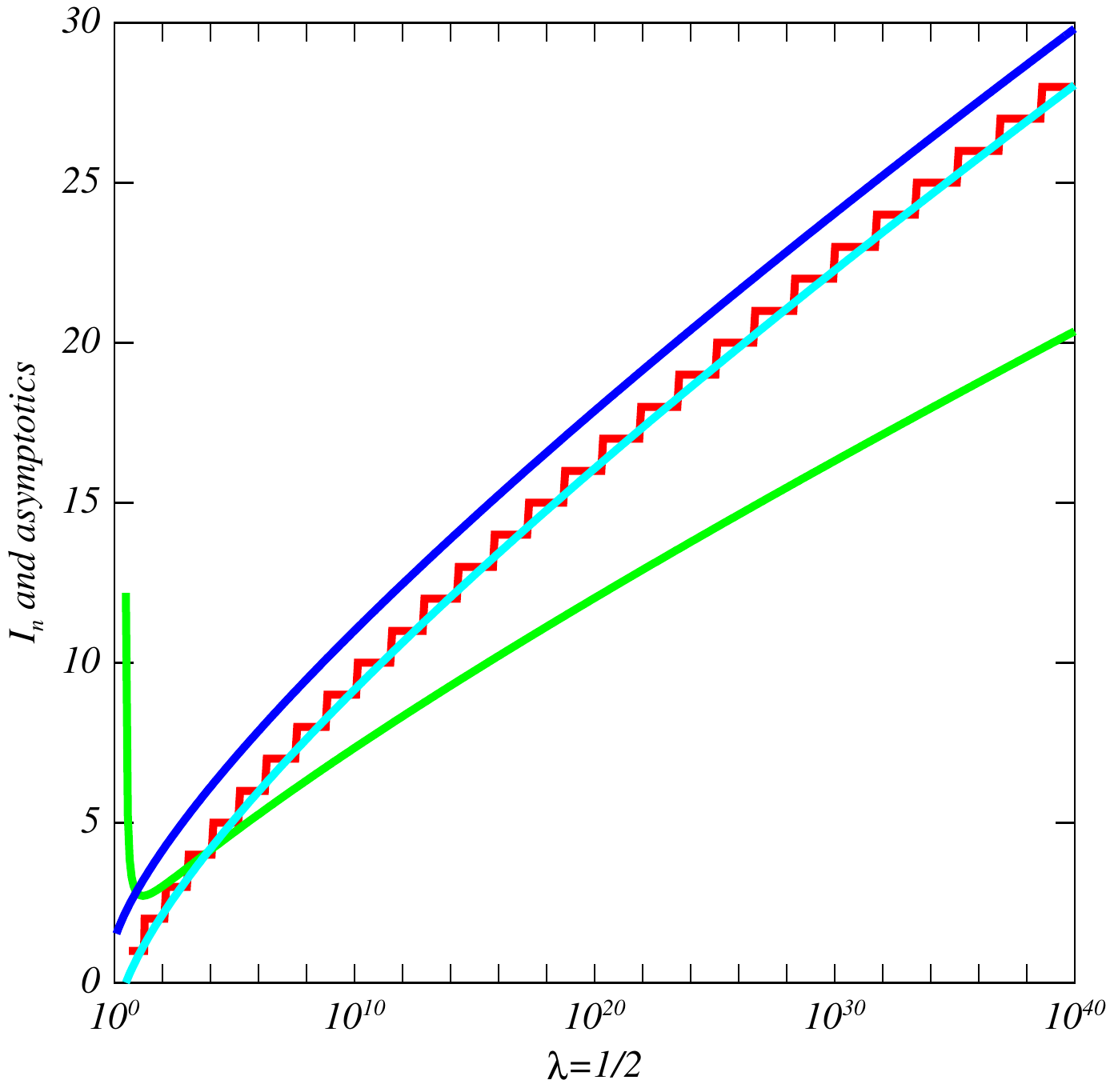}
\includegraphics[width=0.23\linewidth]{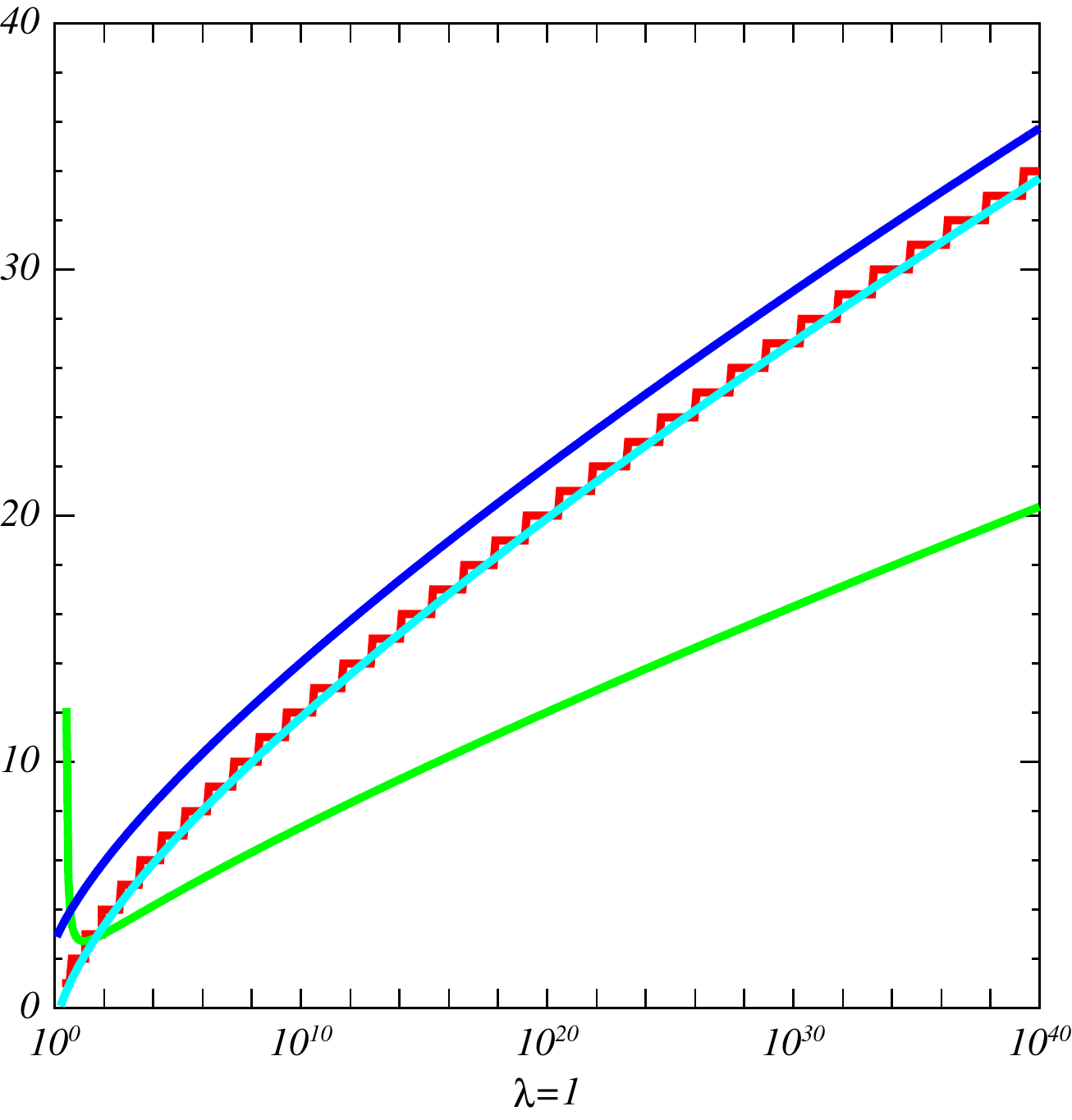}
\includegraphics[width=0.23\linewidth]{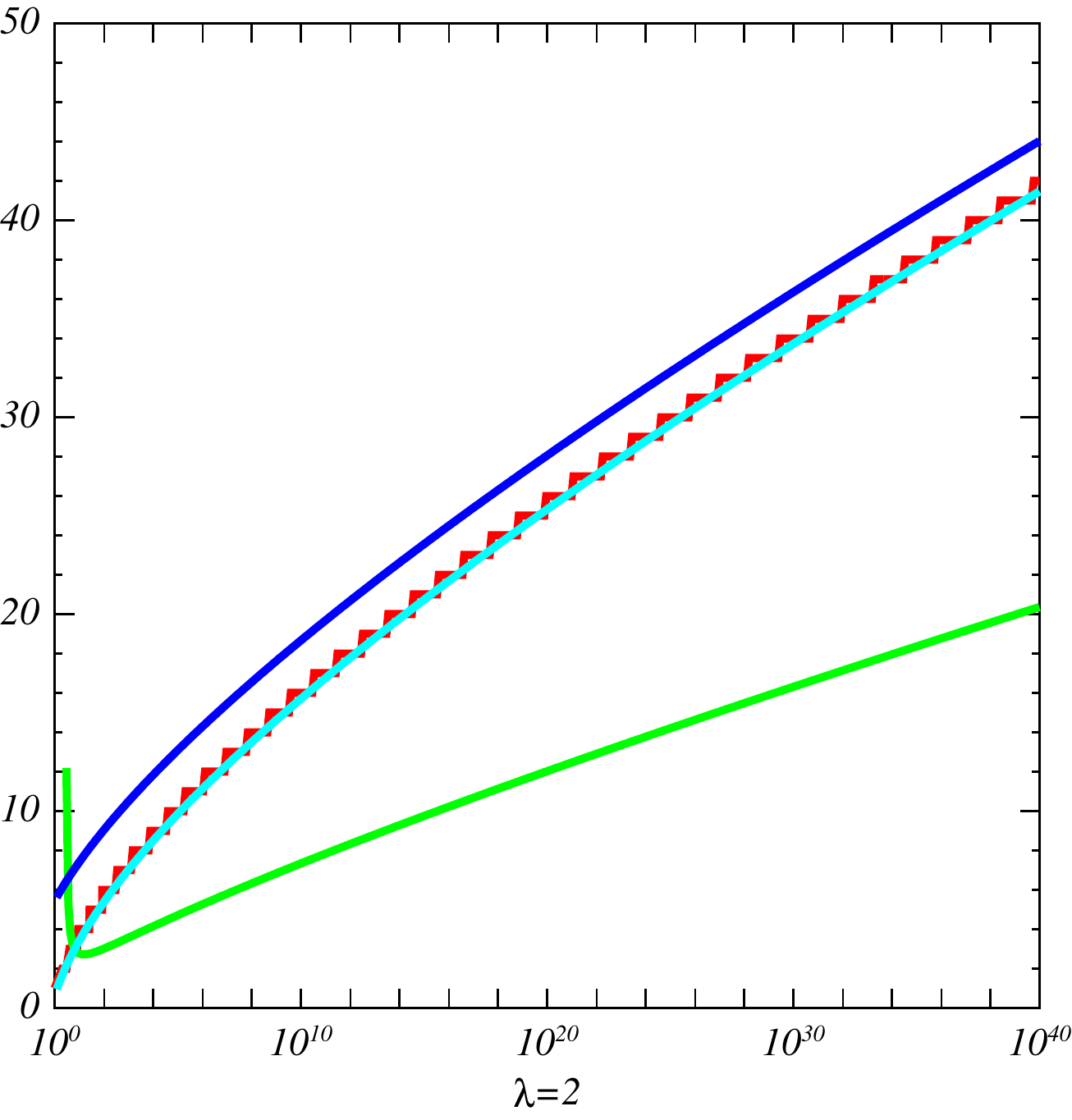}
\includegraphics[width=0.23\linewidth]{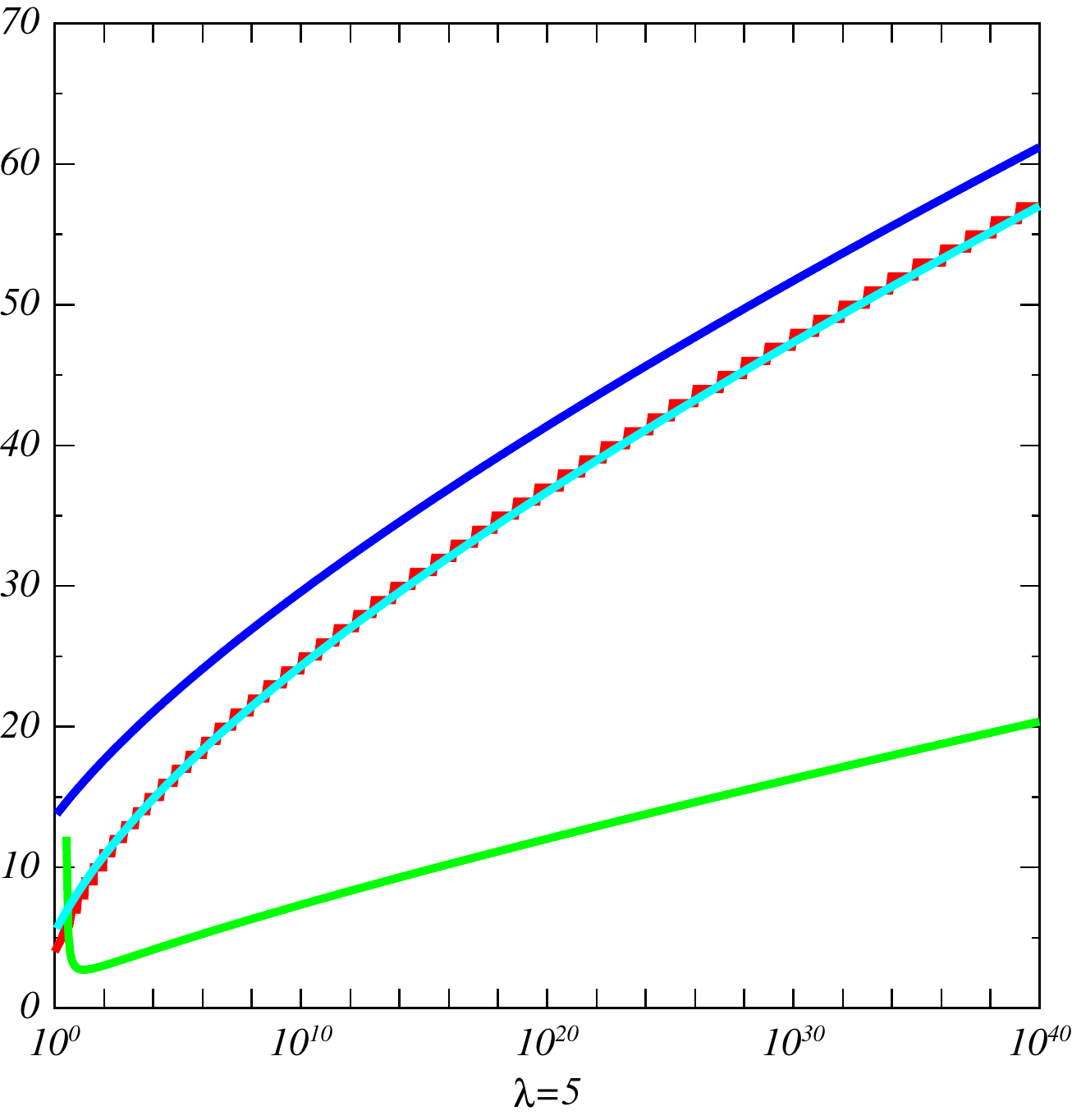}
\end{center}
\caption{Exact values and asymptotics of $I_n$ for $\lambda=1/2, 1, 2, 5$ (left to right) and $n=10^0, \dots, 10^{40}$. The staircase red line (almost hidden by the cyan line) represents the exact mode $I_n$; the other lines represent asymptotic approximations: green for the result of \citet{Kimber1983} (which is independent of $\lambda$), dark blue and cyan for our new results $x_0$ and $x_1$ respectively. The cyan curve always sits between the steps of $I_n$, meaning that $x_1$ has error less than unity.}\label{fig3}
\end{figure}

We have
\begin{equation}
g_\lambda(x)=\exp(-\lambda)\;\lambda^x\;\sum_{i=1}^{\infty}\,\frac{\lambda^i}{\Gamma(x+i+1)}\nonumber
\end{equation}
and we will work with 
\begin{multline}
\log(g_\lambda(x))=-x\log(x)+(1+\log\lambda)x -\frac{3}{2}\log(x)\\+\left(\log\lambda-\lambda-\frac{\log(2\pi)}{2}\right)+\frac{\lambda-13/12}{x}+{\cal O}(x^{-2}).\label{gen}
\end{multline}
A first approximation to the solution of $\log(g_\lambda(x))=-\log n$ 
large and negative is given by keeping only the dominant first two terms in Equation~(\ref{gen}):
\begin{equation}
M_n\sim x_0\equiv\frac{\log n}{W\!\left(\frac{\log n}{\exp(1)\lambda}\right)},\nonumber
\end{equation}
where $W(\cdot)$ is the principal branch of Lambert's $W$ function \citep{W}.  That this is already quite accurate can be seen form the dark blue curves in Figure~\ref{fig3}.  However, we would like to do better; ideally the error should be less than unity so that the mode of the distribution is correctly identified.   A refinement $x_1$ may be generated by making a single Newton correction step; that is, $x_1=x_0-(h(x_0)+\log n)/h'(x_0)$, where $h$ is some approximation to $\log(g_\lambda)$.  For example, by keeping all terms in $\log(g_\lambda(x))$ and $\log(g_\lambda(x))'$  which do not vanish as $n\rightarrow\infty$, we obtain
\begin{equation}
M_n\sim x_1=x_0+\frac{\log\lambda-\lambda-{\log(2\pi)}/{2}-3\log(x_0)/2}{\log(x_0)-\log\lambda}.\nonumber
\end{equation}
This appears to have error less than unity for all values of $n$ and $\lambda$ considered in Figure~\ref{fig3}, and so is probably sufficiently precise for all practical purposes.   If further accuracy is needed, it may be obtained by additional Newton steps.  In any case, both $x_0$ and $x_1$ are considerably more precise than Kimber's approximation.

%
%
%
%


\vfill
\SourceFileFooter

\end{document}